# ANALYTIC CROSSING PROBABILITIES FOR CERTAIN BARRIERS BY BROWNIAN MOTION

By Nabil Kahale

*ESCP-EAP*

We calculate crossing probabilities and one-sided last exit time densities for a class of moving barriers on an interval $[0, T]$ via Schwartz distributions. We derive crossing probabilities and first hitting time densities for another class of barriers on $[0, T]$ by proving a Schwartz distribution version of the method of images. Analytic expressions for crossing probabilities and related densities are given for new explicit and semi-explicit barriers.

**1. Introduction.** Let $W_t$ be a standard Brownian motion and $g$ a continuous function on an interval $[0, T]$. For any real number $u$, consider the first hitting time

$$(1) \qquad \tau(u) = \begin{cases} \inf\{t \in [0, T) : u + W_t = g(t)\}, & \text{if such } t \text{ exists,} \\ T, & \text{otherwise.} \end{cases}$$

Define the boundary-crossing probability with respect to $g$ as $P(\tau(0) < T)$. The calculation of boundary-crossing probabilities or other functionals of the first hitting time density arises in several areas such as statistical testing [7, 23], the valuation of barrier options [9, 14, 24] and default modeling [2, 11]. Closed-form formulae for boundary-crossing probabilities or first time hitting densities are known for few barriers other than the well-known case of linear boundaries. For quadratic boundaries, analytic relations [10, 25] between the first hitting time density and Airy functions have been established. First hitting times for square-root boundaries have been studied in [1, 7, 26]. Analytic expressions involving hypergeometric functions are given in [17] for densities of the first hitting times $\inf\{t \geq 0 : W_t \geq a + b\sqrt{c+t}\}$ and $\inf\{t \geq 0 : |W_t| \geq b\sqrt{c+t}\}$. Using the method of images [4, 5], crossing probabilities and first hitting time densities have been calculated for a class of moving barriers via $\sigma$-finite measures. The method of images is studied in detail in [15], Chapter I. Using the well-known result that $tW_{1/t}$ is a Brownian motion, a correspondence has been established [15], Chapter I, between the method of images and the so-called method of weighted likelihood functions [22, 23] (see also [12], Section 4.3). A relation between the densities of the last exit time and the first hitting time for the time-reversed boundary is shown in [25].

This paper provides crossing probabilities for new moving barriers. Our approach is inspired from the reflection principle (see, e.g., [21], page 105) that shows that a constant barrier-crossing probability is twice the probability that $W_T$ is beyond the barrier. The reflection principle can be easily generalized as follows. If $g(0) > 0$ and $A$ is a Borel subset of $[g(T), +\infty)$ such that the probability that a Brownian motion starting at $g(t)$ at time $t$ belongs to $A$ at time $T$ is a constant $c$ independent of $t$, then $P(\tau(0) < T) = P(W_T \in A)/c$. Theorem 2.5 gives a further generalization of the reflection principle by showing that the crossing probability for certain moving barriers can be retrieved from the density of $W_T$ via Schwartz distributions. The probabilistic basis for the proof of Theorem 2.5 is closely related to Fortet's equation ([8], page 217),

$$\frac{1}{\sqrt{T}} N'\left(\frac{u-v}{\sqrt{T}}\right) = E\left(\frac{1}{\sqrt{T-\tau(u)}} N'\left(\frac{g(\tau(u))-v}{\sqrt{T-\tau(u)}}\right)\right)$$

valid for $u \leq g(0)$ and $v > g(T)$, where $N(x) = (2\pi)^{-1/2} \int_{-\infty}^x \exp(-u^2/2)\,du$. Fortet's equation has been used to establish integral equations [18] and numerical approximations (see [2] and references therein) for the first hitting time density. To our knowledge, though, it has not been previously applied to establish closed-form formulae for curved boundary crossing probabilities. Section 3 gives an analytic expression for the one-sided last exit time density for the barriers considered in Theorem 2.5 under an additional smoothness condition. In Section 4, we show through time-inversion that results in Sections 2 and 3 yield hitting probabilities and first crossing-time densities for certain boundaries on $[T, +\infty)$. Section 5 generalizes our results to two-sided boundaries. Section 6 establishes a Schwartz distribution version of the method of images by building upon our proof techniques. Finally, Section 7 contains concluding remarks. The classical and the Schwartz distribution versions of the method of images allow the calculation of the crossing probability and of the first hitting time density for a class of barriers which is, to our knowledge, unrelated to the one obtained via Theorem 2.5. Several examples illustrate our results. Except for one-sided and double-sided linear barriers, where crossing probabilities and related densities are known or can easily be derived through classical tools, the results in our examples are new, to the best of our knowledge. In particular, we give simple explicit formulae for boundary-crossing probabilities and the one-sided last exit time densities when $g(t) = a + b\sqrt{T-t}$ and, under suitable conditions on $a$ and $b$, when $g(t) = a - \sqrt{(T-t)\ln(\frac{b}{T-t})}$.

**2. Crossing probabilities.** Define the function $H$ from $\mathbb{R}^2$ to $\mathbb{R}^+$ by

(2) $$H(u,v) = P(\tau(u) < T \wedge u + W_T > v).$$

If $\mu$ is a Schwartz distribution or a $\sigma$-finite measure and $u \in \mathbb{R}$, let

$$U(u,t;\mu) = \begin{cases} \int N\left(\frac{u-v}{\sqrt{T-t}}\right) d\mu(v), & \text{if } 0 \leq t < T, \\ 0, & \text{otherwise.} \end{cases}$$

A technical difficulty often encountered in this paper is the interchanging of Schwartz distributions and integration in the measure-theoretic sense, which is valid if certain conditions are met, as shown in Lemma 2.1, but may lead to erroneous conclusions if they are not, as shown in Example 3.

LEMMA 2.1. *Let $m_1$ be a finite measure on $\mathbb{R}$, $m_2$ a Schwartz distribution on $\mathbb{R}$ with compact support $K$, and $\phi(t,v)$ a measurable function on $\mathbb{R}^2$. Assume there exists an open set $V$ containing $K$ and a sequence $\alpha_k$ such that, for any $k \geq 0$, $\frac{\partial^k}{\partial v^k}\phi(t,v)$ exists and is bounded in absolute value by $\alpha_k$ for $(t,v) \in \mathbb{R} \times V$. Then:*

1. *The function $v \mapsto \int \phi(t,v)\,dm_1(t)$ is $C^\infty$ on $V$.*
2. *The function $t \mapsto \int \phi(t,v)\,dm_2(v)$ is bounded and measurable on $\mathbb{R}$.*
3.

(3) $$\int \left(\int \phi(t,v)\,dm_2(v)\right) dm_1(t) = \int \left(\int \phi(t,v)\,dm_1(t)\right) dm_2(v).$$

PROOF. Using the Lebesgue dominated convergence theorem, it can be shown by induction on $k$ that $\frac{\partial^k}{\partial v^k}\int \phi(t,v)\,dm_1(t)$ exists for $v \in V$ and equals $\int \frac{\partial^k}{\partial v^k}\phi(t,v)\,dm_1(t)$. Hence the first assertion.

Since $m_2$ is of finite order and $|\frac{\partial^k}{\partial v^k}\phi(t,v)| \leq \alpha_k$ on $\mathbb{R} \times V$, the function $t \mapsto \int \phi(t,v)\,dm_2(v)$ is bounded on $\mathbb{R}$. To show it is measurable, consider a sequence $\psi_n$ of $C^\infty$ functions with compact support that converges, in the Schwartz distribution sense, to $m_2$. The sequence $\psi_n$ can be chosen so that there exists a compact $K' \subset V$ containing the support of $\psi_n$ for all $n$. Thus, for any $t$, $\int \phi(t,v)\psi_n(v)\,dv \to \int \phi(t,v)\,dm_2(v)$ as $n \to \infty$. Since the function $\phi(t,v)\psi_n(v)$ is bounded and measurable on $\mathbb{R}^2$, the function $t \mapsto \int \phi(t,v)\psi_n(v)\,dv$ is measurable on $\mathbb{R}$ for any $n$. Thus, the function $t \mapsto \int \phi(t,v)\,dm_2(v)$ is measurable on $\mathbb{R}$ as limit of measurable functions.

We now show the third assertion. The measure-theoretic version of Fubini's theorem shows that

$$\int \left( \int \phi(t,v)\psi_n(v)\,dv \right) dm_1(t) = \int \left( \int \phi(t,v)\,dm_1(t) \right) \psi_n(v)\,dv. \tag{4}$$

By the first assertion, the right-hand side of (4) goes to the right-hand side of (3) as $n$ goes to infinity. We prove the third assertion by showing that the left-hand side of (4) goes to the left-hand side of (3) as $n$ goes to infinity.

We say that a set $E$ of $C^\infty$ functions on $\mathbb{R}$ is *bounded* if the supports of all elements of $E$ are contained in the same compact and if there exists a sequence $\beta_k$ such that $|\xi^{(k)}(v)| \leq \beta_k$ for $\xi \in E$, $v \in \mathbb{R}$ and $k \geq 0$. It is known [20], Theorem V.26, that if a sequence of Schwartz distributions converges to a Schwartz distribution, it converges uniformly on any bounded set. Let $\eta$ be a $C^\infty$ function on $\mathbb{R}$ with a compact support contained in $V$ and such that $\eta(v) = 1$ on an open set containing $K \cup K'$. Since the set $E = \{\phi(t,.)\eta, t \in \mathbb{R}\}$ is a bounded set, $\int \phi(t,v)\psi_n(v)\,dv$ converges uniformly to $\int \phi(t,v)\,dm_2(v)$ as $n \to \infty$. As a consequence, the left-hand side of (4) converges to the left-hand side of (3) as $n \to \infty$. $\square$

PROPOSITION 2.2. *For $0 \leq t < T$, the function $N(\frac{v-g(t)}{\sqrt{T-t}})$ is $C^\infty$ with respect to $v$. For any $k \geq 0$ and any closed subset $S$ of $\mathbb{R} - \{g(T)\}$, $\frac{\partial^k}{\partial v^k} N(\frac{v-g(t)}{\sqrt{T-t}})$ is bounded for $(t,v) \in [0,T) \times S$.*

PROOF. For $0 \leq t < T$, $\frac{\partial^k}{\partial v^k} N(\frac{v-g(t)}{\sqrt{T-t}}) = (T-t)^{-k/2} N^{(k)}(\frac{v-g(t)}{\sqrt{T-t}})$, where $N^{(k)}$ is the $k$th derivative of the normal function. Let $\varepsilon > 0$ be the minimum distance between $S$ and $g(T)$. Since $g$ is continuous at $T$, there exists $s \in [0,T)$ such that $|g(t) - g(T)| < \varepsilon/2$ for $s \leq t \leq T$. For $v \in S$ and $0 \leq t < s$, $\frac{\partial^k}{\partial v^k} N(\frac{v-g(t)}{\sqrt{T-t}})$ is bounded by a constant independent of $v$ and $t$ since the function $N^{(k)}$ is bounded. On the other hand, for $v \in S$ and $s \leq t < T$, $|v - g(t)| \geq \varepsilon/2$. Since the function $x^k N^{(k)}(x)$ is bounded, $|\frac{\partial^k}{\partial v^k} N(\frac{v-g(t)}{\sqrt{T-t}})| \leq (\frac{2|v-g(t)|}{\varepsilon\sqrt{T-t}})^k |N^{(k)}(\frac{v-g(t)}{\sqrt{T-t}})|$ is also bounded by a constant independent of $v$ and $t$ for $v \in S$ and $s \leq t < T$. $\square$

PROPOSITION 2.3. *If $u \leq g(0)$ and $v \geq g(T)$, then $H(u,v) = N(\frac{u-v}{\sqrt{T}})$.*

PROOF. If $u + W_T > g(T)$, then $\tau(u) < T$ a.s. since $W_t$ is a.s. a continuous function of $t$. Thus, for $v \geq g(T)$,

$$P(\tau(u) < T \wedge u + W_T > v) = P(u + W_T > v)$$
$$= N\left(\frac{u-v}{\sqrt{T}}\right). \qquad \square$$

LEMMA 2.4. 1. *For fixed $u$, the function $H(u,\cdot)$ is $C^\infty$ on $\mathbb{R} - \{g(T)\}$.*
2. *If $\nu$ is a Schwartz distribution with a compact support contained in $\mathbb{R} - \{g(T)\}$ then, for any real number $u$,*

$$\int H(u,v)\,d\nu(v) = E(U(g(\tau(u)), \tau(u); \nu)). \tag{5}$$

PROOF. Let $m_1$ be the image measure of the random variable $\tau(u)$ and $\phi(t,v) = N(\frac{g(t)-v}{\sqrt{T-t}})$ if $0 \le t < T$, and $\phi(t,v) = 0$ otherwise. We can write

$$\begin{aligned}
H(u,v) &= P(\tau(u) < T \wedge u + W_T > v) \\
&= E(1_{\tau(u)<T} P(u + W_T > v | \mathcal{F}_{\tau(u)})) \\
&= E(1_{\tau(u)<T} P(W_T - W_{\tau(u)} > v - g(\tau(u)) | \mathcal{F}_{\tau(u)})) \\
&= E(\phi(\tau(u), v)) \\
&= \int \phi(t,v)\,dm_1(t).
\end{aligned} \tag{6}$$

The second equality follows from the tower law and the fourth one from the strong Markov property. It follows from (6), Proposition 2.2 and Lemma 2.1 that $H(u,v)$ is $C^\infty$ with respect to $v$ on any open subset $V$ of $\mathbb{R}$ such that $g(T) \notin \overline{V}$. Hence, the function $H(u,\cdot)$ is $C^\infty$ on $\mathbb{R} - \{g(T)\}$.

We show the second assertion:

$$\begin{aligned}
\int H(u,v)\,d\nu(v) &= \int \left( \int \phi(t,v)\,dm_1(t) \right) d\nu(v) \\
&= \int \left( \int \phi(t,v)\,d\nu(v) \right) dm_1(t) \\
&= \int U(g(t), t; \nu)\,dm_1(t) \\
&= E(U(g(\tau(u)), \tau(u); \nu)).
\end{aligned}$$

The first equality follows from (6) and the second one from (3). □

THEOREM 2.5. *Let $\mu_n$ be a sequence of Schwartz distributions with compact supports contained in $[g(T), +\infty)$ such that, for $0 \le t < T$, the sequence $U(g(t), t; \mu_n) \to 1$ as $n \to \infty$ and, for some constant $c$, $|U(w, t; \mu_n)| \le c$ for all $n$, $0 \le t < T$ and $w \le g(t)$. Then, for $u \le g(0)$, $U(u, 0; \mu_n) \to P(\tau(u) < T)$ as $n \to \infty$.*

PROOF. Given $n$ and $\varepsilon > 0$, consider the Schwartz distribution $\nu_{n,\varepsilon}$ such that $\int \psi(v)\,d\nu_{n,\varepsilon}(v) = \int \psi(v+\varepsilon)\,d\nu_n(v)$ for any $C^\infty$ function $\psi$ with compact support on $\mathbb{R}$. The support of $\nu_{n,\varepsilon}$ is contained is $(g(T), +\infty)$ and so, by (5),

$$\int H(u,v)\,d\nu_{n,\varepsilon}(v) = E(U(g(\tau(u)), \tau(u); \nu_{n,\varepsilon})). \tag{7}$$

But it follows from Proposition 2.3 that the left-hand side of (7) is equal to $U(u, 0; \nu_{n,\varepsilon})$. On the other hand, $U(x, t; \nu_{n,\varepsilon}) = U(x - \varepsilon, t; \nu_n)$. By letting $\varepsilon$ go to 0 and using the Lebesgue dominated convergence theorem, (7) becomes

$$U(u, 0; \nu_n) = E(U(g(\tau(u)), \tau(u); \nu_n)).$$

By applying once again the Lebesgue dominated convergence theorem, it follows that $U(u, 0; \nu_n) \to P(\tau(u) < T)$ as $n \to \infty$. $\square$

COROLLARY 2.6. *Let $\mu$ be a Schwartz distribution with a compact support included in $[g(T), +\infty)$ such that, for some constant $c$ and $0 \leq t < T$, $U(g(t), t; \mu) = 1$ and $|U(w, t; \mu)| \leq c$ if $w \leq g(t)$. Then $P(\tau(u) < T) = U(u, 0; \mu)$ for $u \leq g(0)$.*

COROLLARY 2.7. *Let $\mu$ be a $\sigma$-finite measure on $[g(T), +\infty)$ such that $U(g(t), t; \mu) = 1$ for $0 \leq t < T$. Then $P(\tau(u) < T) = U(u, 0; \mu)$ for $u \leq g(0)$.*

REMARK 2.8. *Let $0 \leq t < T$. A simple time-change argument shows that, for $u \leq g(t)$, the probability that a Brownian motion starting at $u$ at time $t$ hits the barrier $g$ between $t$ and $T$ is equal to $\lim_{n \to \infty} U(u, t; \mu_n)$.*

**3. One-sided last exit time densities.** Define

$$\lambda = \sup\{t \in (0, T] : W_t = g(t)\}$$

and

$$\sigma = \sup\{t \in (0, T] : W_t \geq g(t)\}.$$

Observe that

$$P(\sigma \leq t) = P(\lambda \leq t \wedge W_T \leq g(T)).$$

Denote by

$$N_2(x, y; \rho) = \int_{-\infty}^{x} \int_{-\infty}^{y} \frac{1}{2\pi \sqrt{1 - \rho^2}} \exp\left(-\frac{u^2 + v^2 - 2\rho uv}{2(1 - \rho^2)}\right) du\, dv$$

the cumulative bivariate normal probability distribution function.

THEOREM 3.1. *Under the assumptions of Theorem 2.5, for $0 < t \leq T$,*

$$(8) \quad P(\sigma \leq t) = N\left(\frac{g(t)}{\sqrt{t}}\right) - \lim_{n \to \infty} \int N_2\left(\frac{g(t)}{\sqrt{t}}, -\frac{v}{\sqrt{T}}; -\sqrt{\frac{t}{T}}\right) d\mu_n(v).$$

*Furthermore, if $g$ is $C^1$ on an open subset $I$ of $[0, T]$ and the sequence $U(g(t), t; \mu'_n)$ is uniformly bounded and converges to a continuous function of $t$ as $t \in I$ and $n \to \infty$, then*

$$(9) \quad \frac{d}{dt} P(\sigma \leq t) = \frac{1}{2\sqrt{t}} N'\left(\frac{g(t)}{\sqrt{t}}\right) \lim_{n \to \infty} U(g(t), t; \mu'_n)$$

*for $t \in I$.*

PROOF. Let $t \in (0, T)$. By Remark 2.8,

$$P(W_t \leq g(t) \wedge \lambda > t) = \int_{-\infty}^{g(t)} \frac{1}{\sqrt{t}} N'\left(\frac{u}{\sqrt{t}}\right) \lim_{n \to \infty} U(u, t; \mu_n) \, du$$

$$= \lim_{n \to \infty} \int_{-\infty}^{g(t)} \frac{1}{\sqrt{t}} N'\left(\frac{u}{\sqrt{t}}\right) U(u, t; \mu_n) \, du$$

$$= \lim_{n \to \infty} \int \left( \int_{-\infty}^{g(t)} \frac{1}{\sqrt{t}} N'\left(\frac{u}{\sqrt{t}}\right) N\left(\frac{u-v}{\sqrt{T-t}}\right) du \right) d\mu_n(v)$$

$$= \lim_{n \to \infty} \int N_2\left(\frac{g(t)}{\sqrt{t}}, -\frac{v}{\sqrt{T}}; -\sqrt{\frac{t}{T}}\right) d\mu_n(v).$$

The second equation follows from the Lebesgue dominated convergence theorem and the third one from Lemma 2.1 with $m_1(du) = \frac{du}{1+u^2}$, $m_2 = \mu_n$, $V = \mathbb{R}$ and $\phi(u, v) = (1+u^2)\frac{1}{\sqrt{t}} N'(\frac{u}{\sqrt{t}}) N(\frac{u-v}{\sqrt{T-t}}) 1_{u < g(t)}$. Since

$$P(\sigma \leq t) = P(\lambda \leq t \wedge W_t \leq g(t))$$
$$= P(W_t \leq g(t)) - P(W_t \leq g(t) \wedge \lambda > t),$$

we conclude that (8) holds. On the other hand, since

$$\frac{\partial N_2}{\partial x}(x, y; \rho) = N'(x) N\left(\frac{y - \rho x}{\sqrt{1 - \rho^2}}\right)$$

and [19]

$$\frac{\partial N_2}{\partial \rho}(x, y; \rho) = \frac{1}{\sqrt{1 - \rho^2}} N'(x) N'\left(\frac{y - \rho x}{\sqrt{1 - \rho^2}}\right),$$

(10)
$$\frac{\partial}{\partial t} N_2\left(\frac{g(t)}{\sqrt{t}}, -\frac{v}{\sqrt{T}}; -\sqrt{\frac{t}{T}}\right)$$
$$= N'\left(\frac{g(t)}{\sqrt{t}}\right) \left( N\left(\frac{g(t) - v}{\sqrt{T-t}}\right) \frac{\partial}{\partial t}\left(\frac{g(t)}{\sqrt{t}}\right) \right.$$
$$\left. - \frac{1}{2\sqrt{t(T-t)}} N'\left(\frac{g(t) - v}{\sqrt{T-t}}\right) \right).$$

Let $J_n(t) = \int N_2(\frac{g(t)}{\sqrt{t}}, -\frac{v}{\sqrt{T}}; -\sqrt{\frac{t}{T}}) d\mu_n(v)$. If $[s, t] \subset I$ with $s < t$,

(11) $\quad J_n(t) - J_n(s) = \int_s^t \left( \int \frac{\partial}{\partial \theta} N_2\left(\frac{g(\theta)}{\sqrt{\theta}}, -\frac{v}{\sqrt{T}}; -\sqrt{\frac{\theta}{T}}\right) d\mu_n(v) \right) d\theta.$

Equation (11) can be derived by applying Lemma 2.1 with $m_1(d\theta) = 1_{s \leq \theta \leq t} \, d\theta$, $m_2 = \mu_n$, $\phi(\theta, v) = 1_{s \leq \theta \leq t} \frac{\partial}{\partial \theta} N_2(\frac{g(\theta)}{\sqrt{\theta}}, -\frac{v}{\sqrt{T}}; -\sqrt{\frac{\theta}{T}})$ and $V = \mathbb{R}$. Combining

(10) and (11) yields

$$J_n(t) - J_n(s)$$
$$= \int_s^t N'\left(\frac{g(\theta)}{\sqrt{\theta}}\right)\left(U(g(\theta),\theta;\mu_n)\frac{\partial}{\partial \theta}\left(\frac{g(\theta)}{\sqrt{\theta}}\right) - \frac{1}{2\sqrt{\theta}}U(g(\theta),\theta;\mu_n')\right)d\theta$$

and so, by the Lebesgue dominated convergence theorem,

$$\lim_{n\to\infty} J_n(t) - J_n(s)$$
(12)
$$= N\left(\frac{g(t)}{\sqrt{t}}\right) - N\left(\frac{g(s)}{\sqrt{s}}\right)$$
$$- \int_s^t N'\left(\frac{g(\theta)}{\sqrt{\theta}}\right) \lim_{n\to\infty} \frac{1}{2\sqrt{\theta}}U(g(\theta),\theta;\mu_n')\,d\theta.$$

Equation (9) follows by combining (8) and (12). □

EXAMPLE 1. Let $\mu_n(v) = 2\delta_{a+bT} + 2be^{2b(v-a-bT)}1_{\{a+bT<v<n\}}$ for $n \geq 0$ and $g(t) = a + bt$. Then

$$U(u,t;\mu_n) \to N\left(\frac{u-g(T)}{\sqrt{T-t}}\right) + e^{2b(u-g(t))}N\left(\frac{u+g(T)-2g(t)}{\sqrt{T-t}}\right)$$

as $n$ goes to infinity. By Theorem 2.5,

(13) $$P(\tau(0) < T) = N\left(\frac{-a-bT}{\sqrt{T}}\right) + e^{-2ab}N\left(\frac{bT-a}{\sqrt{T}}\right)$$

if $a \geq 0$. By Theorem 3.1, for $0 < t < T$,

$$\frac{d}{dt}P(\sigma \leq t) = \frac{1}{\sqrt{t}}\left(\frac{N'(b\sqrt{T-t})}{\sqrt{T-t}} + bN(b\sqrt{T-t})\right)N'\left(\frac{g(t)}{\sqrt{t}}\right).$$

By symmetry,

$$\frac{d}{dt}P(\lambda \leq t)$$
(14)
$$= \frac{1}{\sqrt{t}}\left(\frac{2N'(b\sqrt{T-t})}{\sqrt{T-t}} + b(N(b\sqrt{T-t}) - N(-b\sqrt{T-t}))\right)$$
$$\times N'\left(\frac{g(t)}{\sqrt{t}}\right)$$

for $0 < t < T$. Note that (13) agrees with the well-known formula [27], Corollary 7.2.2, for linear barriers crossing. When $a = b = 0$, (14) is identical to the well-known first arcsine law [21], page 112.

EXAMPLE 2. $\mu = \frac{1}{N(b)}\delta_a$ and $g(t) = a + b\sqrt{T-t}$. By Corollary 2.7

$$P(\tau(0) < T) = N\left(\frac{-a}{\sqrt{T}}\right) \Big/ N(b), \tag{15}$$

if $a + b\sqrt{T} \geq 0$. By Theorem 3.1, for $0 < t < T$ and any values of $a$, $b$,

$$P(\sigma \leq t) = N\left(\frac{g(t)}{\sqrt{t}}\right) - \frac{1}{N(b)} N_2\left(\frac{g(t)}{\sqrt{t}}, -\frac{a}{\sqrt{T}}; -\sqrt{\frac{t}{T}}\right)$$

and

$$\frac{d}{dt} P(\sigma \leq t) = \frac{N'(b)}{2N(b)\sqrt{t(T-t)}} N'\left(\frac{g(t)}{\sqrt{t}}\right).$$

By symmetry,

$$P(\lambda > t)$$
$$= \frac{1}{N(b)} N_2\left(\frac{g(t)}{\sqrt{t}}, -\frac{a}{\sqrt{T}}; -\sqrt{\frac{t}{T}}\right) + \frac{1}{N(-b)} N_2\left(\frac{-g(t)}{\sqrt{t}}, \frac{a}{\sqrt{T}}; -\sqrt{\frac{t}{T}}\right)$$

and

$$\frac{d}{dt} P(\lambda \leq t) = \frac{N'(b)}{2N(b)N(-b)} \frac{1}{\sqrt{t(T-t)}} N'\left(\frac{g(t)}{\sqrt{t}}\right) \tag{16}$$

for $0 < t < T$. Equation (15) agrees with the well-known formula [21], page 105, for constant barriers crossing when $b = 0$. Here again, (16) is identical to the first arcsine law when $a = b = 0$.

EXAMPLE 3. $\mu = \sqrt{2\pi b}\delta'_a$, where $b \geq T$. The conditions of Corollary 2.6 hold for $g(t) = a - \sqrt{(T-t)\ln(\frac{b}{T-t})}$. Thus $P(\tau(0) < T) = \sqrt{b/T}\exp(-\frac{a^2}{2T})$ if $a \geq \sqrt{T\ln(b/T)}$. By Theorem 3.1, for any value of $a$ and $0 < t < T$,

$$P(\sigma \leq t) = N\left(\frac{g(t)}{\sqrt{t}}\right) - \sqrt{\frac{b}{T}} \exp\left(-\frac{a^2}{2T}\right) N\left(\frac{g(t) - at/T}{\sqrt{t(1-t/T)}}\right)$$

and

$$\frac{d}{dt} P(\sigma \leq t) = \frac{1}{2}\sqrt{\frac{1}{t(T-t)} \ln \frac{b}{T-t}} N'\left(\frac{g(t)}{\sqrt{t}}\right).$$

Observe that the function $g^*(t) = a + \sqrt{(T-t)\ln(\frac{b}{T-t})}$ does satisfy the equation $U(g^*(t), t; \mu) = 1$ but $U(w, t; \mu)$ is not bounded for $0 \leq t < T$ and $w \leq g^*(t)$. Thus, neither the conditions nor the conclusion of Corollary 2.6 hold for $g^*$.

EXAMPLE 4. More generally, consider the Hermite polynomial $H_n$ defined by

$$\frac{d^n}{dx^n}\exp(-x^2/2) = (-1)^n H_n(x)\exp(-x^2/2).$$

Thus $H_0(x) = 1$, $H_1(x) = x$, $H_2(x) = x^2 - 1$ and $H_{n+1}(x) = xH_n(x) - H'_n(x)$ for $n \geq 0$. Moreover, $H_n$ has $n$ distinct zeros, and two successive zeros of $H_n$ enclose a single zero of $H_{n-1}$. For $n \geq 1$, let $z_n$ be the largest zero of $H_n$ and let $\mu = \sqrt{2\pi}b\delta_a^{(n)}$, with $b \geq \exp(z_n^2/2)T^{n/2}/H_{n-1}(z_n)$. For $0 \leq t < T$, let $x(t)$ be the largest solution to the equation $H_{n-1}(x)\exp(-x^2/2) = (T-t)^{n/2}/b$ and $g(t) = a - \sqrt{T-t}x(t)$, with $g(T) = a$. By Corollary 2.6, if $0 \leq g(0)$, which is equivalent to the conditions

$$a \geq z_n\sqrt{T} \quad \text{and} \quad bT^{-n/2}H_{n-1}\left(\frac{a}{\sqrt{T}}\right)\exp\left(-\frac{a^2}{2T}\right) \leq 1,$$

then

$$P(\tau(0) < T) = bT^{-n/2}H_{n-1}\left(\frac{a}{\sqrt{T}}\right)\exp\left(-\frac{a^2}{2T}\right).$$

For any value of $a$ and $0 < t < T$,

$$\frac{d}{dt}P(\sigma \leq t) = \frac{H_n(x(t))}{2H_{n-1}(x(t))\sqrt{t(T-t)}}N'\left(\frac{g(t)}{\sqrt{t}}\right).$$

When $n = 2$, $g(t) = a - \sqrt{T-t}x(t)$, where $x(t)$ is the largest solution to the equation $x\exp(-x^2/2) = (T-t)/b$. The function $g(t)$ is, up to a translation and time-reversion, identical to the upper boundary of a symmetric band used in [13] as a confidence region for Brownian motion. Whether our techniques could be used to compute analytically the crossing probability for the two-sided band, which has been approximated numerically in [13], is a question that deserves further investigation. Note that when $n = 2$ or $n = 3$, the function $g$ can be expressed in terms of the Lambert function defined as the multivalued inverse of the function $w \mapsto we^w$. This has already been observed in [13] when $n = 2$. Efficient algorithms for computing the Lambert function can be found in [3].

**4. Time-inversion.** If $W_t$ is a standard Brownian motion, then

$$\hat{W}_t = \begin{cases} 0, & \text{if } t = 0, \\ tT^{-1}W_{T^2/t}, & \text{if } t > 0, \end{cases}$$

is also ([21], page 21) a standard Brownian motion. By letting $\hat{g}(t) = tT^{-1}g(T^2/t)$ and

$$\hat{\sigma} = \inf\{t \geq T : W_t \geq \hat{g}(t)\},$$

it follows that

(17) $$P(\exists t \geq T : W_t = \hat{g}(t)) = P(\exists t \in (0,T] : W_t = g(t))$$

and, for $t > T$,

(18) $$\frac{d}{dt}P(\hat{\sigma} \leq t) = -\frac{d}{dt}P(\sigma \leq T^2/t)$$

if the latter derivative exists.

EXAMPLE 5. It follows from Example 2, (17) and (18) that, for $\hat{g}(t) = at/T + b\sqrt{t^2/T - t}$,

$$P(\exists t \geq T : W_t = \hat{g}(t)) = N\left(\frac{-a}{\sqrt{T}}\right) \Big/ N(b)$$

if $a + b\sqrt{T} > 0$. For $t > T$ and any values of $a$, $b$,

$$\frac{d}{dt}P(\hat{\sigma} \leq t) = \frac{N'(b)}{2N(b)t\sqrt{t/T - 1}}N'\left(\frac{\hat{g}(t)}{\sqrt{t}}\right).$$

**5. Two-sided boundaries.** Theorems 2.5 and 3.1 can be extended to two-sided boundaries as follows. Let $g_0$ and $g_1$ be two continuous functions on $[0,T]$ such that $g_1(t) < g_0(t)$ for $0 \leq t \leq T$. Define

$$\tau(u) = \inf\{t \in [0,T) : u + W_t = g_0(t) \text{ or } u + W_t = g_1(t)\}$$

and

$$\sigma = \sup\{t \in (0,T] : W_t \geq g_0(t) \text{ or } W_t \leq g_1(t)\}.$$

If $\mu$ is a Schwartz distribution or a $\sigma$-finite measure, $u \in \mathbb{R}$ and $0 \leq t < T$, let $U_0(u,t;\mu) = U(u,t;\mu)$ and

$$U_1(u,t;\mu) = \int N\left(\frac{v-u}{\sqrt{T-t}}\right) d\mu(v).$$

THEOREM 5.1. *Let $\mu_{0,n}$ (resp. $\mu_{1,n}$) be a sequence of Schwartz distributions with compact supports contained in $[g_0(T), +\infty)$ (resp. $(-\infty, g_1(T)]$) such that, for $0 \leq t < T$, $w = g_0(t)$ or $w = g_1(t)$,*

$$\sum_{i=0}^{1} U_i(w,t;\mu_{i,n}) \to 1$$

*as $n \to \infty$ and, for $g_1(t) \leq w \leq g_0(t)$, $i \in \{0,1\}$, all $n$ and some constant $c'$, $|U_i(w,t;\mu_{i,n})| \leq c'$. Then, for $g_1(0) \leq u \leq g_0(0)$,*

$$P(\tau(u) < T) = \lim_{n \to \infty} \sum_{i=0}^{1} U_i(u,t;\mu_{i,n}).$$

*Furthermore, for $0 < t \leq T$,*

$$P(\sigma \leq t)$$
$$= \lim_{n \to \infty} \sum_{i=0}^{1} (-1)^i \bigg( N\bigg(\frac{g_i(t)}{\sqrt{t}}\bigg) - \int N_2\bigg(\frac{g_i(t)}{\sqrt{t}}, -\frac{v}{\sqrt{T}}; -\sqrt{\frac{t}{T}}\bigg) d\mu_{0,n}(v)$$
$$- \int N_2\bigg(\frac{g_i(t)}{\sqrt{t}}, \frac{v}{\sqrt{T}}; \sqrt{\frac{t}{T}}\bigg) d\mu_{1,n}(v) \bigg).$$

*Finally, if $g_i$ is $C^1$ on an open subset $I$ of $[0,T]$, for $i \in \{0,1\}$, and the sequence $U(g_i(t), t; \mu'_{0,n} - \mu'_{1,n})$ is uniformly bounded and converges to a continuous function of $t$ as $t \in I$ and $n \to \infty$, then*

$$\frac{d}{dt} P(\sigma \leq t) = \sum_{i=0}^{1} \frac{(-1)^i}{2\sqrt{t}} N'\bigg(\frac{g_i(t)}{\sqrt{t}}\bigg) \lim_{n \to \infty} U(g_i(t), t; \mu'_{0,n} - \mu'_{1,n})$$

*for $t \in I$.*

EXAMPLE 6. Let $g_0(t) = a$ and $g_1(t) = b$ for $0 \leq t \leq T$, where $b < a$. Define $\mu_{0,n} = 2 \sum_{j=0}^{n} (-1)^j \delta_{a+j(a-b)}$ and $\mu_{1,n} = 2 \sum_{j=0}^{n} (-1)^j \delta_{b+j(b-a)}$. Then

$$U_i(w, t; \mu_{i,n}) = 2 \sum_{j=0}^{n} (-1)^j N\bigg(\frac{(-1)^i(w - g_i(0)) + j(b-a)}{\sqrt{T-t}}\bigg)$$

is positive and upper-bounded by 1 for $i \in \{0,1\}$ and $b \leq w \leq a$. If $w = a$ or $w = b$, then

$$\sum_{i=0}^{1} U_i(w, t; \mu_{i,n}) = 1 + 2(-1)^n N\bigg(\frac{(n+1)(b-a)}{\sqrt{T-t}}\bigg) \to 1$$

as $n \to \infty$. If $b < 0 < a$, Theorem 5.1 yields the well-known formula (see, e.g., [21], Exercise 3.15, page 111, for a direct derivation)

$$P(\tau(0) < T) = 2 \sum_{j=0}^{\infty} (-1)^j N\bigg(\frac{-a + j(b-a)}{\sqrt{T}}\bigg)$$
$$+ 2 \sum_{j=0}^{\infty} (-1)^j N\bigg(\frac{b + j(b-a)}{\sqrt{T}}\bigg).$$

Furthermore, for $0 < t < T$,

$$\frac{d}{dt} P(\sigma \leq t) = \frac{1}{\sqrt{2\pi t(T-t)}} \bigg( N'\bigg(\frac{a}{\sqrt{t}}\bigg) + N'\bigg(\frac{b}{\sqrt{t}}\bigg) \bigg)$$
$$\times \bigg( 1 + 2 \sum_{j=1}^{\infty} (-1)^j \exp\bigg(-\frac{j^2(a-b)^2}{2(T-t)}\bigg) \bigg).$$

EXAMPLE 7. Let $\mu_0 = c^{-1}\delta_a$ and $\mu_1 = c^{-1}\delta_b$, where $b < a$ and $2N(\frac{b-a}{2\sqrt{T}}) < c < 1$. For $0 \leq t \leq T$, let $g_0(t)$ be the largest solution to the equation

$$N\left(\frac{w-a}{\sqrt{T-t}}\right) + N\left(\frac{b-w}{\sqrt{T-t}}\right) = c,$$

and $g_1(t) = a + b - g_0(t)$. If $b < 0 < a$, then, by Theorem 5.1,

$$P(\tau(0) < T) = c^{-1}N\left(\frac{-a}{\sqrt{T}}\right) + c^{-1}N\left(\frac{b}{\sqrt{T}}\right)$$

and

$$\frac{d}{dt}P(\sigma \leq t) = \frac{1}{2c\sqrt{t(T-t)}}\left(N'\left(\frac{g_0(t)}{\sqrt{t}}\right) + N'\left(\frac{g_1(t)}{\sqrt{t}}\right)\right)$$
$$\times \left(N'\left(\frac{g_0(t)-a}{\sqrt{T-t}}\right) - N'\left(\frac{g_0(t)-b}{\sqrt{T-t}}\right)\right).$$

**6. Schwartz distributions and the method of images.** Let $B_t$, $t \geq 0$, be a Brownian motion and $f$ a continuous function on an interval $[0, T']$ such that $f(0) > 0$. Define the first hitting time $\tau$ by

$$\tau = \inf\{t \in [0, T'): B_t = f(t)\}.$$

Fix a real number $T \in (0, T']$ and let $h$ be the function from $\mathbb{R}$ to $\mathbb{R}^+$ defined by

$$h(u) = P(\tau < T \wedge B_T \leq u).$$

For $0 \leq t \leq T$, let $g(t) = f(T-t)$ and $W_t = B_{T-t} - B_T$. Thus, $W_t$ is a Brownian motion [21], page 21, on $[0, T]$ and $W_T = -B_T$. As before, the random variable $\tau(u)$ and the function $H$ are defined via (1) and (2), respectively. The following proposition shows formally that the probability that a Brownian motion starting at $u$ at time 0 and ending at 0 at time $T$ hits $g$ equals the probability that a Brownian motion starting at 0 at time 0 and ending at $u$ at time $T$ hits $f$.

PROPOSITION 6.1. *For $u < g(0)$, $h'(u) = -\frac{\partial H}{\partial v}(u, 0)$.*

PROOF. By Lemma 2.4, $H(u, \cdot)$ has a first derivative with respect to $v$ at $v = 0$ since $g(T) = f(0) > 0$. Similarly, using the relation $h(u) = E(1_{\tau < T}N(\frac{u-f(\tau)}{\sqrt{T-\tau}}))$, it can be shown that $h$ is $C^\infty$ on $\mathbb{R} - \{g(0)\}$. For $\varepsilon > 0$,

(19)
$$P(\tau(u) < T \wedge 0 < u + W_T \leq \varepsilon)$$
$$= P(\tau(u) < T \wedge u - \varepsilon \leq B_T < u)$$
$$\geq P(\tau < T \wedge u - \varepsilon \leq B_T < u).$$

To show the second equation, we observe that if $\tau < T$ and $u - \varepsilon \leq B_T < u$, then $B_t = f(t)$ for some $t \in (0, T)$. Thus, $W_{T-t} + B_T = g(T - t)$ and so $W_{T-t} + u > g(T - t)$. Since $u < g(0)$ and $W$ is a.s. continuous on $[0, T]$, we conclude that, a.s., $W_{t'} + u = g(t')$ for some $t' \in [T - t, T)$. Hence $\tau(u) < T$ a.s. It follows from (19) that $H(u, 0) - H(u, \varepsilon) \geq h(u) - h(u - \varepsilon)$. By letting $\varepsilon \to 0$, we conclude that $h'(u) \leq -\frac{\partial H}{\partial v}(u, 0)$.

Similarly,

$$\begin{aligned}P(\tau(u) < T \wedge \ -\varepsilon < u + W_T \leq 0) \\ = P(\tau(u) < T \wedge \ u \leq B_T < u + \varepsilon) \\ \leq P(\tau < T \wedge \ u \leq B_T < u + \varepsilon).\end{aligned} \tag{20}$$

We show the second equation by noting that if $\tau(u) < T$ and $u \leq B_T < u + \varepsilon$, then $u + W_t = g(t)$ for some $t \in (0, T)$. Thus, $u + B_{T-t} - B_T = f(T - t)$ and $B_{T-t} \geq f(T - t)$. Since $f(0) > 0$, we conclude that, a.s., $B_{t'} = f(t')$ for some $t' \in (0, T - t]$, and so $\tau < T$. It follows from (20) that $H(u, -\varepsilon) - H(u, 0) \leq h(u + \varepsilon) - h(u)$, and so $h'(u) \geq -\frac{\partial H}{\partial v}(u, 0)$. Thus $h'(u) = -\frac{\partial H}{\partial v}(u, 0)$, as desired. □

THEOREM 6.2. *Let $\mu$ be a Schwartz distribution with compact support contained in the interval $(f(0), +\infty)$. Assume that $\int \exp(\frac{2f(t)v - v^2}{2t}) d\mu(v) = 1$ for $t \in (0, T']$. Then*

$$P(\tau < T) = N\left(\frac{-f(T)}{\sqrt{T}}\right) + U(f(T), 0; \mu). \tag{21}$$

*If $f$ is $C^1$ on an open subset $I$ of $[0, T']$, then, for $T \in I$,*

$$\frac{\partial}{\partial T} P(\tau < T) = \frac{1}{2T^{3/2}} \int v N'\left(\frac{f(T) - v}{\sqrt{T}}\right) d\mu(v). \tag{22}$$

PROOF. Consider the Schwartz distribution $\nu = \delta'_0 - \mu'$. For $0 \leq t < T$,

$$\begin{aligned}U(g(t), t; \nu) &= U(f(T - t), t; \delta'_0 - \mu') \\ &= \frac{1}{\sqrt{T-t}}\left(N'\left(\frac{f(T-t)}{\sqrt{T-t}}\right) - \int N'\left(\frac{f(T-t) - v}{\sqrt{T-t}}\right) d\mu(v)\right) \\ &= 0.\end{aligned}$$

It follows from Lemma 2.4 that $\int H(u, v) d\nu(v) = 0$ for any real number $u$ and so, by Proposition 2.3,

$$\begin{aligned}\frac{\partial H}{\partial v}(u, 0) &= \int \frac{\partial H}{\partial v}(u, v) d\mu(v) \\ &= -\int \frac{1}{\sqrt{T}} N'\left(\frac{u - v}{\sqrt{T}}\right) d\mu(v) \\ &= -U(u, 0; \mu').\end{aligned} \tag{23}$$

Using Proposition 6.1, we conclude that $h'(u) = U(u, 0; \mu')$ for $u < f(T)$. Since $h$ is a continuous function of $u$ and is $C^1$ on $(-\infty, f(T))$, it follows that, for $u < f(T)$,

$$
\begin{aligned}
h(u) &= \int_{-\infty}^{u} h'(x)\, dx \\
&= \int_{-\infty}^{u} \left( \int \frac{1}{\sqrt{T}} N'\left(\frac{x-v}{\sqrt{T}}\right) d\mu(v) \right) dx \\
&= \int \left( \int_{-\infty}^{u} \frac{1}{\sqrt{T}} N'\left(\frac{x-v}{\sqrt{T}}\right) dx \right) d\mu(v) \\
&= U(u, 0; \mu).
\end{aligned}
\tag{24}
$$

The third equation follows from Lemma 2.1 with $m_1(dx) = \frac{dx}{1+x^2}$, $m_2 = \mu$, $\phi(x, v) = (1+x^2) N'(\frac{x-v}{\sqrt{T}}) 1_{x<u}$ and $V = (f(0), +\infty)$.

We infer from (24) that $P(\tau < T \wedge B_T \leq f(T)) = U(f(T), 0; \mu)$. On the other hand, since the events $\tau < T \wedge B_T > f(T)$ and $B_T > f(T)$ coincide up to a negligible set, $P(\tau < T \wedge B_T > f(T)) = N(-f(T)/\sqrt{T})$. This concludes the proof of (21).

Assume now that $f$ is $C^1$ on an open subset $I$ of $[0, T']$ containing $T$. Then

$$
\begin{aligned}
\frac{\partial}{\partial T} P(\tau < T) &= \int \frac{\partial}{\partial T}\left(\frac{f(T)-v}{\sqrt{T}}\right) N'\left(\frac{f(T)-v}{\sqrt{T}}\right) d\mu(v) \\
&\quad - \frac{\partial}{\partial T}\left(\frac{f(T)}{\sqrt{T}}\right) N'\left(\frac{f(T)}{\sqrt{T}}\right) \\
&= \int \frac{\partial}{\partial T}\left(\frac{-v}{\sqrt{T}}\right) N'\left(\frac{f(T)-v}{\sqrt{T}}\right) d\mu(v) \\
&= \frac{1}{2T^{3/2}} \int v N'\left(\frac{f(T)-v}{\sqrt{T}}\right) d\mu(v).
\end{aligned}
$$

The first equality can be justified by fixing a real number $s$ such that $[s, T] \subset I$ and applying Lemma 2.1 with $m_1(dt) = 1_{s \leq t \leq T}\, dt$, $m_2 = \mu$, $V = \mathbb{R}$ and $\phi(t, v) = 1_{s \leq t \leq T} \frac{\partial}{\partial t} N(\frac{f(t)-v}{\sqrt{t}})$. □

EXAMPLE 8. $\mu = b\delta'_a$, where $a > 0$ and $b \geq a \exp(1 - \frac{a^2}{2T'})$. Let $f(0) = a/2$ and $f(t) = a + ty(t)/a$ for $0 < t \leq T'$, where $y(t)$ is the smallest solution to the equation $y(t)\exp(y(t)) = -a\exp(-\frac{a^2}{2t})/b$. The function $f$ can be easily expressed in terms of the Lambert function. By Theorem 6.2, it follows after some calculations that, for $0 < T < T'$,

$$
P(\tau < T) = N\left(\frac{-f(T)}{\sqrt{T}}\right) + \frac{\sqrt{T}}{a - f(T)} N'\left(\frac{f(T)}{\sqrt{T}}\right)
$$

and

$$\frac{d}{dT}P(\tau < T) = \frac{1}{2T^{3/2}}\bigg(a - \frac{T}{a - f(T)}\bigg)N'\bigg(\frac{f(T)}{\sqrt{T}}\bigg).$$

**7. Concluding remarks.** We conclude with the following remarks:

1. Our results can be adapted to other Markov processes, including discrete-time Markov chains.
2. The condition $|U(w,t;\mu_n)| \leq c$ for all $n$, $0 \leq t < T$ and $w \leq g(t)$ in Theorem 2.5 can be relaxed as follows: $|U(w,t;\mu_n)| \leq c$ for all $n$, $s \leq t < T$ and $w \leq g(t)$, for some fixed $s \in [0,T)$.
3. Another approach to prove Theorem 2.5 is to show that $U(u,t;\mu)$ satisfies the partial differential equation

$$\frac{\partial U}{\partial t} + \frac{1}{2}\frac{\partial^2 U}{\partial u^2} = 0$$

and to use the martingale properties of $U(W_t, t; \mu)$.
4. The Schwartz distributions in our examples were essentially found by inspection. The questions of establishing necessary and sufficient conditions for the existence of a sequence of Schwartz distributions satisfying the conditions of Theorem 2.5 for a given function $g$, and determining such a sequence if it exists, are left for future research. A numerical method for calculating $\sigma$-finite measures associated to a given boundary via the classical method of images has been used in [6, 16] to approximate the first hitting time density in cases where the exact density is not known or difficult to compute.

**Acknowledgments.** The author would like to thank Monique Jeanblanc, Eric Séré and an anonymous referee for constructive comments.

ESCP-EAP  
79 avenue de la République  
75543 Paris Cedex 11  
France  
E-mail: nkahale@escp-eap.net  
URL: http://nkahale.free.fr